\documentclass[12pt]{amsart}

\usepackage{amsthm}
\usepackage{amsmath}
\usepackage{amssymb}
\usepackage{eucal}
\usepackage{amsfonts}
\usepackage{amsmath}
\usepackage[Lenny]{fncychap}
\usepackage{enumerate}
\usepackage{amssymb}
\usepackage[english]{babel}
\usepackage[breaklinks=true]{hyperref}
\usepackage[all]{xy}
\usepackage{fancybox}
\usepackage{latexsym,mathrsfs,longtable,lscape,dsfont,yfonts}
\usepackage{srcltx}
\usepackage{soul}

\ifx\pdfpagewidth\undefined 
\usepackage[T1]{fontenc}
\else 
\usepackage[OT1]{fontenc} 
\fi 
\usepackage{amsfonts,amssymb,latexsym,longtable,amsmath}

\newtheorem{thm}{Theorem}[section]
\newcommand{\bthm}{\begin{thm}} \newcommand{\ethm}{\end{thm}}
\newtheorem{prop}[thm]{Proposition}
\newcommand{\bprp}{\begin{prop}} \newcommand{\eprp}{\end{prop}}
\newtheorem{fact}[thm]{Fact}
\newcommand{\bfct}{\begin{fact}} \newcommand{\efct}{\end{fact}}
\newtheorem{prob}[thm]{Problem}
\newcommand{\bprb}{\begin{prob}} \newcommand{\eprb}{\end{prob}}
\newtheorem{quest}[thm]{Question}
\newcommand{\bqtn}{\begin{quest}} \newcommand{\eqtn}{\end{quest}}
\newtheorem{lem}[thm]{Lemma}
\newcommand{\blem}{\begin{lem}} \newcommand{\elem}{\end{lem}}
\newtheorem{claim}[thm]{Claim}
\newcommand{\bclm}{\begin{claim}} \newcommand{\eclm}{\end{claim}}
\newtheorem{cor}[thm]{Corollary}
\newcommand{\bcor}{\begin{cor}} \newcommand{\ecor}{\end{cor}}
\newtheorem{conj}[thm]{Conjecture}
\newcommand{\bcnj}{\begin{conj}} \newcommand{\ecnj}{\end{conj}}
\theoremstyle{definition}
\newtheorem{defn}[thm]{Definition}
\newcommand{\bdfn}{\begin{defn}} \newcommand{\edfn}{\end{defn}}
\newtheorem{spec}[thm]{Specializing}
\newcommand{\bspc}{\begin{spec}} \newcommand{\espc}{\end{spec}}
\theoremstyle{remark}
\newtheorem{rem}[thm]{Remark}
\newcommand{\brem}{\begin{rem}} \newcommand{\erem}{\end{rem}}
\newtheorem{cnv}[thm]{Convention}
\newcommand{\bcnv}{\begin{cnv}} \newcommand{\ecnv}{\end{cnv}}
\newtheorem{exam}[thm]{Example}
\newcommand{\bexm}{\begin{exam}} \newcommand{\eexm}{\end{exam}}
\newcommand{\bpf}{\begin{proof}} \newcommand{\epf}{\end{proof}}

\newtheorem{thmy}{\textbf{Theorem}}

\newcommand{\C}{\mathbb C}
\newcommand{\Z}{\mathbb Z}

\newcommand{\Q}{\mathbb Q}

\renewcommand{\phi}{\varphi}
\renewcommand{\theta}{\vartheta}

\newcommand{\gd}{{\delta}}

\newcommand{\w}{{\rm w}}

\newcommand{\U}{\mathbb{U}}

\newcommand{\mkp}{\medskip}

\def\defi{\buildrel\rm def \over=}

\begin{document}

\title[The weak compactification of locally compact groups]{The weak compactification of locally compact groups}

\author[M. Ferrer]{Mar\'ia V. Ferrer}
\address{Universitat Jaume I, Instituto de Matem\'aticas de Castell\'on,
Campus de Riu Sec, 12071 Castell\'{o}n, Spain.}
\email{mferrer@mat.uji.es}

\author[S. Hern\'andez]{Salvador Hern\'andez}
\address{Universitat Jaume I, Departamento de Matem\'{a}ticas,
Campus de Riu Sec, 12071 Castell\'{o}n, Spain.}
\email{hernande@mat.uji.es}

\thanks{ Research Partially supported by the Spanish Ministerio de Econom\'{i}a y Competitividad,
grant: MTM/PID2019-106529GB-I00 (AEI/FEDER, EU)}
\vspace{1cm}

\begin{abstract}
We further investigate the weak topology generated by the irreducible unitary representations of a group $G$.
A deep result due to Ernest \cite{Ernest1971} and Hughes \cite{Hughes1973} asserts that every weakly
compact subset of a locally compact (LC) group $G$ is compact in the LC-topology, generalizing thereby a
previous result of Glicksberg \cite{glicks1962} for abelian locally compact (LCA) groups. Here, we first survey some recent findings
on the weak topology and establish some new results about the preservation of several compact-like properties when going from the weak topology
to the original topology of LC groups. Among others, we deal with the preservation
of countably compactness, pseudocompactness and functional boundedness.

\end{abstract}


\dedicatory{Respectfully dedicated to Professor \emph{Hans-Peter Kunzi} }

\date{\today}

\maketitle \setlength{\baselineskip}{24pt}

\section{Introduction}

 In contrast to what happens with abelian groups, where the Bohr compactification and Bohr topology display many nice features, the Bohr compactification presents important shortcomings when applied to non-commutative groups.
For example, it may happen that the Bohr compactification of a locally compact group
becomes trivial, what makes it useless in order to study the structure of those non-abelian groups.
On the other hand, it is known, as a consequence of the celebrated Gel'fand and Ra\v{\i}kov
Theorem, that the set of all irreducible unitary representations of a locally compact group $G$ contains all the information necessary
to recover the topological and algebraic structure of the group (see \cite{EnockSchwartz1992}).
Therefore, it seems appropriate to consider the weak topology generated by the irreducible unitary representations of a group $G$
as the genuine weak topology of general not ne\-ces\-sa\-rily LC group $G$ (as a matter of fact, they coincide
for abelian topological groups, by the Schur's Lemma). This is what we have done in \cite{FerHerTar_JFA} and, here, we further develop this approach, which was initiated by Hughes in \cite{Hughes1973}.

There is a plethora of results that concern the weak topology of abelian topological groups and, although there still are interesting open questions in this setting, we have plenty of information about the weak topology of abelian LC groups. The literature in this regard is vast, so we only mention \cite{Comfort2001} and the references therein.
On the other hand, there also are some crucial findings
that have been already established for general LC groups, which are less known so far but form the basis
for an in-depth study of weak topologies in the setting of not necessarily abelian groups.
It is pertinent to mention here the task developed by Ernest \cite{Ernest1971} and Hughes \cite{Hughes1973}, who proved that every weakly compact subset
of a LC group $G$ is compact in the LC-topology, generalizing thereby a previous result of
Glicksberg \cite{glicks1962} for abelian locally compact (LCA) groups (see also \cite{remustrig99} where the Bohr topology of 
Moore groups is investigated). The main goal of this paper is twofold. First,
we survey some recent results on the weak topology of LC groups. On the other hand, we further develop this research line
by studying the preservation of several compact-like properties when going from the weak topology to the original topology
of LC groups. Among others, we deal with the preservation of countably compactness, pseudocompactness and functional boundedness.
We now collect some definitions and basic facts that will be used along the paper.
\section{Basic facts}
\subsection{The Bohr topology}
The Bohr compactification of a topological group is a well known notion that has been widely treated in the setting of topological groups.
Nevertheless, we remind here its most basic features for the reader's sake.
With every (not necessarily abelian) topological group $ G $ there is associated a compact Hausdorff group
$ bG $, the so-called \emph{Bohr compactification} of $ G $, and a continuous homomorphism
$ b $ of $ G $ onto a dense subgroup of $ bG $ such that $ bG $ is characterized by the following universal property:
given any continuous homomorphism $ h $ of $ G $ into a compact group $ K $,
there is always a continuous homomorphism $ \bar {h} $ of $ bG $ into $ K $ such that
$ h = \bar {h} \circ b $ (see \cite [V \S 4]{Heyer1970}, where a detailed study on $bG$ and their properties is given).
The \emph{Bohr topology} of a topological group $G$ is the one that inherits as a subgroup of $bG$.

In the sequel, following a terminology introduce by Trigos-Arrieta \cite{Trigos-Arrieta1991},
if $G$ is a topological group, we denote by $G^+$ the same algebraic group but equipped with the Bohr topology.
However, when the group $G$ is discrete, we will use the symbolism $G^\sharp$ used by van Douwen in \cite{glicks1962}.

As far as we know, the first main result about the Bohr topology is due to Glicksberg \cite{glicks1962}.

\bthm[Glicksberg, 1962]\label{Glicksberg} Let $(G,\tau )$ be a LCA group. It
holds that $A\subseteq G$ is $\tau $-compact if and only if is $\tau ^{+}$-compact.
\ethm
\mkp

However, the systematic study of the Bohr topology was started by van Douwen (loc. cit.)
in the setting of discrete abelian groups.

\bthm[van Douwen, 1990]\label{vdfirst}
Every $A\subseteq G$ contains a subset $D$ with $|D|=|A|$ that is
relatively discrete and $C$-embedded in $G^\sharp$ as well as
$C^\ast$-embedded in $bG$.
\ethm

In addition to the previous main result, van Douwen studied the Bohr topology of a discrete abelian group in depth.
It is remarkable that, except for the standing abelian
hypothesis, his proofs of results concerning {\large$\sharp$}-groups
made no use whatsoever of specific algebraic properties. This
probably led him to ask  whether two groups $G_1$ and $G_2$ with the
same cardinality should have $G_1^{\sharp}$ and $G_2^{\sharp}$
homeomorphic. Some years later Kunen \cite{Kunen1998} and,
independently, { Dikranjan} and { Watson} \cite{DikranWat2001},  gave
examples of {\em torsion} groups with the same cardinality having
nonhomeomorphic $\sharp$-spaces. Still, much remains unknown, even
among groups of countable cardinality. One can get an idea of how
involved the situation is by taking into account, see
\cite{cht3}, that $\Q^{\sharp}$ and
$((\Q/\Z)\times\Z)^{\sharp}$ are homeomorphic.

Van Douwen's work on the Bohr topology of discrete abelian groups was continued by Trigos-Arrieta, for locally compact abelian (LCA) groups,
in his doctoral dissertation \cite{Trigos-Arrieta1991}. There, the author introduces the following notion:

\bdfn[Trigos-Arrieta, 1991]\label{Respect}
Let $G$ be a topological group and let $\mathcal P$ be a topological property.
We say that $G$ \emph{respects} $\mathcal P$ if for any subset $F$ of $G$ the following holds:
the subspace $F$ of $G$ has $\mathcal P$ if and only if the subspace $F$ of $G^+$ has $\mathcal P$.
\edfn
\mkp

For example, since $G=G^+$ for $G$ compact, it is obvious that compact groups respects all topological properties.
In \cite{Trigos-Arrieta1991}, Trigos-Arrieta proves that LCA groups respect most compact-like properties:
pseudocompactness, functional boundedness, and other topological properties: Lindel\"ofness and connectedness.

Further properties concerning the Bohr
topology of a LCA group can be found in
\cite{comfherntrig,cht3,gal-her-fm1999,Hernandez1998}

\subsection{The weak topology of LC groups}
Given a locally compact group $(G,\tau)$, we denote by $Irr(G)$ the set of all continuous
unitary irreducible representations $\sigma$ defined on $G$.
That is, continuous in the sense that each matrix coefficient function
$g\mapsto \langle\sigma(g)u,v\rangle$ is a continuous map of $G$ into the complex plane.
Thus, fixed $\sigma\in Irr(G)$, if $\mathcal{H}^{\sigma}$ denotes the Hilbert space associated to $\sigma$,
we equip the unitary group $\U(\mathcal H^\sigma)$ with the weak (equivalently, strong) operator topology.
For two elements $\pi$ and $\sigma$ of $Irr(G)$, we write $\pi\sim\sigma$
to denote the relation of unitary equivalence and we denote by $\widehat G$ the \emph{dual object} of $G$,
which is defined as the set of equivalence classes in ($Irr(G)/{\sim}$).
We refer to \cite{Dixmier1964,Folland,Bekka2008} for all undefined
notions concerning the unitary representations of locally compact groups.

Adopting, the terminology introduced by Ernest in \cite{Ernest1971}, set $\mathcal{H}_n\defi \C^n$ for $n=1,2, \ldots$;
and $\mathcal H_{0}\defi l^2(\Z)$. The symbol $Irr^C_n(G)$ will denote the set of irreducible unitary representations of $G$
on $\mathcal{H}_n$, where it is assumed that every set $Irr^C_n(G)$ is equipped with the compact open topology.
Finally, define $Irr^C(G)=\bigsqcup\limits_{n\geq 0} Irr^C_n(G)$ (the disjoint topological sum).

We denote by $G^{\rm w}=(G,{\rm w}(G,Irr(G))$ 
the group $G$ equipped with the weak (group) topology generated by $Irr(G)$. 
Since equivalent representations define the same topology, we have $G^{\rm w}=(G,{\rm w}(G,\widehat G))$. That is, the \emph{weak topology}
is the initial topology on $G$ defined by the dual object. Moreover, in case $G$ is a separable, metric, locally compact group,
then every irreducible unitary representation acts on a separable Hilbert space and, as a consequence, is unitary equivalent
to a member of $Irr^C(G)$. Thus $G^{\rm w}=(G,{\rm w}(G,Irr^C(G)))$ for separable, metric, locally compact groups.
We will make use of this fact in order to avoid
the proliferation of isometries (see \cite{Dixmier1964}).
In case the group $G$ is abelian,
the dual object $\widehat G$ is a group, which is called \emph{dual group},
and the weak topology of $G$ reduces to the weak topology generated by
all continuous homomorphisms of $G$ into the unit circle $\mathbb T$. That is,
the weak topology coincides with the {Bohr topology} of $G$.

\subsection{The weak compactification of LC groups}
\bdfn\label{I(G)}
We denote by $P(G)$ the set of continuous positive definite functions on $(G,\tau)$.
If $\sigma\in Irr(G)$ and $v\in \mathcal{H}^{\sigma}$, 
then the positive definite function:
$$\varphi:g\mapsto \langle\sigma (g)v,v\rangle\text{,   }  g\in G$$
is called \textit{pure}, and the family of all such functions is denoted by $I(G)$.
When $G$ is abelian, the set $I(G)$ coincides with the dual group $\widehat G$ of the group $G$.
\edfn
\mkp

The proof of the lemma below is straightforward.

\blem\label{top_debil_I}
Let $G$ be a locally compact group. Then
$G^{{\rm w}} = (G, {\rm w} (G,I(G)))$.
\elem

\mkp

\bdfn\label{compactification wG}
Let $G$ be a locally compact group and consider the following natural embedding:
\begin{equation*}
\w :G\hookrightarrow \prod\limits_{\varphi\in I(G)}\overline{\varphi(G)}\ \
\hbox{with}\ \ \w (g)=(\varphi(g))_{\varphi\in I(G)}
\end{equation*}

\mkp

The \emph{weak compactification} $\w G$ of $G$ as the pair $(\w G,\w)$,
where $\w G\defi \overline{\w (G)}$.

This compactification has been previously considered in \cite{Cheng2011,Cheng2013} using different techniques. Also Akemann and Walter \cite{Akemann1972}
extended Pontryagin duality to non-abelian locally compact groups using the family of pure positive definite functions.
Again, in case $G$ is abelian, the compactifican $(\w G,\w)$ coincides with $bG$, the
Bohr compactification of  $G$.

A better known compactification of a locally group $G$ which is closely related to $\w G$ is
defined as follows (cf. \cite{Eymard1964,Spronk2013}): let $\overline{\mathrm{B(G)}}^{\|\cdot\|_\infty}$ denote
the commutative $C^*$-algebra consisting of the uniform closure of the Fourier-Stieltjes algebra
of $G$. Here, the Fourier-Stieltjes algebra is defined as the matrix coefficients of
the unitary representations of $G$. Following \cite{Mayer1997} we call the spectrum $eG$ of
$\overline{\mathrm{B(G)}}^{\|\cdot \|_\infty}$ the
\emph{Eberlein compactification} of $G$.
Since the Eberlein compactification $eG$ is defined using the family of all continuous positive definite functions, it follows
that $\w G$ is a factor of $eG$ and, as a consequence,  inherits most of its properties. In particular,
$\w G$ is a compact involutive semitopological semigroup. 
\edfn
\mkp

We now recall some known results about unitary representations of locally compact groups that are needed in the proof of our main result in this section.
One main point is the decomposition of unitary representations by direct integrals of irreducible unitary representations.
This was established by Mautner \cite{Mautner1950} following the ideas introduced by von Neuman in \cite{vonNeumann1949}.

\bthm[F. I. Mautner, \cite{Mautner1950}]\label{Th_Mautner}
For any representation $(\sigma,\mathcal H_\sigma)$ of a separable locally compact group $G$, there is a measure space
$(R,\mathcal R,r)$, a family $\{\sigma[p]\}$ of irreducible representations of $G$, which are associated
to each $p\in R$, and an isometry $U$ of $\mathcal H_\sigma$ such that
$$U\sigma U^{-1}=\int_R \sigma[p]d_rp.$$
\ethm

\brem\label{Re_integral decomposition}
The proof of the above theorem given by Mautner assumes that the representation space $\mathcal H_\sigma$ is separable
but, subsequently, Segal \cite{Segal1951} removed this constraint. Furthermore, it is easily seen that
we can assume that $\sigma[p]$ belongs to $Irr^C(G)$ locally almost everywhere in the theorem above (cf. \cite{Ikeshoji1979}).
\erem

A remarkable consequence of Theorem \ref{Th_Mautner}
is the following corollary about positive definite functions.

\bcor\label{Co_Mautner}
Every Haar-measurable positive definite function $\varphi$ on a separable locally compact group $G$ can be
expressed for all $g\in G$ outside a certain set of Haar-measure zero in the form
$$\varphi(g)= \int_R \varphi_p(g)d_rp,$$

\noindent where $\varphi_p$ is a pure positive definite functions on $G$ for all $p\in R$.
\ecor
\mkp

The following proposition is contained in the proof of Lemma 3.2 of Bichteler \cite[pp. 586-587]{Bichteler1969}

\bprp\label{Pro_referee}
Let $G$ be a locally compact group. If $H$ is an open subgroup of $G$, then
each continuous irreducible representation of $H$
is the restriction of a continuous irreducible representation of $G$.
\eprp

\section{Locally compact groups respect compactness}
In this section we prove an old result by Ernest \cite{Ernest1971} (cf. \cite{Ikeshoji1979}) and \cite{Hughes1973},
asserting that the weak topology of locally compact groups respects compactness, that is it holds that
every weakly compact subset is compact for the original locally compact group topology.
In fact, Ernest first proved that, for separable metric locally compact groups, every weakly convergent sequence
is convergent for the locally compact topology and this result was subsequently extended by Hughes by proving that
the weak topology of locally compact groups respects compactness.

Unfortunately, even though the formulation of Hugues' result quoted above can be found in \cite{Hughes1973}, its full proof only appears
in his Doctoral dissertation but has never been published later.
Therefore, we have decided to include it here for the reader's sake. Our proof is complete since it contains both Ernest's and Hughes' results.
We first need a further definition and several previous lemmas.

\bdfn
Let $U$ be an open neighbourhood of the identity of a topological group $G$. We say that a sequence
$\lbrace g_n\rbrace_{n<\omega}$ is \textit{$U$-discrete} if $g_nU\cap g_mU=\emptyset$ for all $n\neq m\in\omega$.
\edfn

\blem\emph{(J. Ernest, 1971)}\label{Lem_sequences}
Let $(G,\tau)$ be a separable metric locally compact group. Then every convergent sequence $\{g_{n}\}_{n<\omega}$ in $G^{\rm w}$ is
also $\tau$-convergent in $G$.
\elem
\bpf
Remark that there is no loss of generality in assuming that $\{g_{n}\}_{n<\omega}$ converges to $e_G$ in $G^{\rm w}$.
We must verify that $\lbrace g_{n}\rbrace_{n<\omega }$ is $\tau$-convergent to $e_G$.

In order to do so, take a continuous positive definite function $\psi\in P(G)$. Since $G$ is separable, by Corollary \ref{Co_Mautner},
there is a measure space $(R,\mathcal R,r)$, a family  $\{\psi_p\}$ of pure positive definite functions on $G$, which are associated
to each $p\in R$, such that
$$\psi(g)= \int_R \psi_p(g)d_rp\ \hbox{for all}\ g\in G.$$

Therefore

$$\psi(g_{n})= \int_R \psi_p(g_{n})d_rp\ \hbox{for all}\ n<\omega.$$

Now, for each $n<\omega$, consider the map $f_n$ on $R$ by $f_n(p)\defi \psi_p(g_{n})$.
Then $f_n$ is integrable on $R$ and, since  $\lbrace g_{n}\rbrace_{n<\omega }$ is weakly convergent to $e_G$,
it follows that $\{f_n(p)\}$ converges to $\psi_p(e_G)$ for all $p\in R$. Furthermore, if
$\psi_p(g)=<\sigma_p(g)v_p,v_p>$ for some $\sigma_p\in Irr(G)$ and $v_p\in\mathcal H_{\sigma_p}$,
it follows that $$|f_n(p)|=|\psi_p(g_{n})|=|<\sigma_p(g_{n})v_p,v_p>|\leq \|v_p\|^2.$$
Thus  defining $f$ on $R$ as the pointwise limit of $\{f_n\}$,
we are in position to apply Lebesgue's dominated convergence theorem in order to obtain that
$$\psi(e_G)= \int_R \psi_p(e_{G})d_rp=\int_R f(p)d_rp=\lim\limits_{n\rightarrow \infty} \int_R \psi_p(g_{n})d_rp=\lim\limits_{n\rightarrow \infty} \psi(g_{n}).$$

In other words, the sequence $\{\psi(g_{n})\}$ converges to $\psi(e_G)$ for all $\psi\in P(G)$.
Hence $\{g_{n}\}$ $\tau$-converges to $e_G$ in $G$ and we are done. 
(see \cite{GelfandRaikov} or \cite[Prop. 3.33]{Folland}).
\epf

\bthm\emph{(J.R. Hughes, 1972)}\label{compact_3}
\emph{Let $(G,\tau)$ be a locally compact group. Then $(G,\tau)$ and $G^{\rm w}$ contain the same compact subsets.}
\ethm
\bpf
Let $B$ be a weakly compact subset of $G$. Remark that if $B$ were $\tau$-precompact, since $B$ is $\tau$-closed in $G$, it would follow
that $B$ is $\tau$-compact in $G$.

Thus, reasoning by contradiction, we assume that $B$ is not $\tau$-precompact in $G$.
Then there exists an open, symmetric and relatively compact neighbourhood of the identity $U$ in $G$
such that $B$ contains a $U$-discrete sequence $\lbrace g_n\rbrace_{n<\omega}$.

Consider the subgroup $H\defi < \overline U\cup\lbrace g_n\rbrace_{n<\omega} >$, which $\sigma$-compact and open in $G$.
By Kakutani-Kodaira's theorem, there exists a normal, compact subgroup $K$ of $H$ such that $K\subseteq U$ and $H/K$ is metrizable, and consequently separable (Polish).
Let $p:H\rightarrow H/K$ be the quotient homomorphism and let $\overline p :\w H\rightarrow \w\, H/K$ denote its canonical extension
to the weak compactifications. By Proposition \ref{Pro_referee}, we have that
$\overline H^{\, \w G}$ is canonically homeomorphic to $\w H$. Therefore $\{g_{n}\}_{n<\omega}$ weakly converges to the neutral element in $H$.
Hence $\{p(g_{n})\}_{n<\omega}$ weakly converges to the neutral element in $H/K$, which is a separable, metrizable, LC group.
Thus, by Lemma \ref{Lem_sequences}, the sequence $\{p(g_{n})\}_{n<\omega}$ $\tau/K$-converges to the neutral element in $H/K$.
Then by a theorem of Varopoulos \cite{Varopoulos1964}, the sequence $\lbrace p(g_n)\rbrace_{n<\omega}$ can be lifted to
a sequence $\lbrace x_n\rbrace_{n<\omega}\subseteq H$ converging to some point $x_0\in H$. This entails that $ x_n^{-1}g_n\in K$ for all $n\in\omega$.
Thus the sequence $\lbrace g_n\rbrace_{n<\omega}$ would be contained in the compact subset
$(\lbrace x_n\rbrace_{n<\omega}\cup\lbrace x_0\rbrace) K$, which is a contradiction since $\lbrace g_n\rbrace_{n<\omega}$ was supposed
to be $U$-discrete. This contradiction completes the proof.
\epf

\subsection{Weakly Cauchy sequences}

In some special cases, Hughes' theorem implies the convergence
of weakly Cauchy sequences. Indeed, let us denote by $inv(\w G)$ the group of invertible elements of $\w G$.
It is known (see [38, Proposition II.4.6.(i)]) that every maximal subgroup of a compact
semitopological semigroup is a topological group that is complete with respect to the two-sided uniformity.
In particular, this applies to $inv(\w G)$, which is a complete (for the two-sided uniformity) topological group.

\bprp\label{compact_4}
Let $(G,\tau)$ be a locally compact group and suppose that $\{g_{n}\}_{n<\omega}$ is a Cauchy sequence in $G^{\rm w}$.
If $\overline{\{g_{n}\}_{n<\omega}}^{\, \w G}\subseteq inv(\w G)$, then $\{g_{n}\}_{n<\omega}$ is $\tau$-convergent in $G$.
\eprp
\bpf
Assume that $\{g_{n}\}_{n<\omega}$ is a Cauchy sequence in $G^{\rm w}$. First, we verify that the sequence is
a precompact subset of $(G,\tau)$.

Indeed, we have that $\{g_{n}\}_{n<\omega}$ converges to some element $p\in inv(\w G)$.
If $\{g_{n}\}_{n<\omega}$ were not precompact in $(G,\tau)$, there would be a neighbourhood of the
neutral element $U$ and a subsequence $\{g_{n(m)}\}_{m<\omega}$ such that
$g_{n(m)}^{-1}\cdot g_{n(l)}\notin U$ for each $m,l<\omega$ with $m\neq l$.
On the other hand, since $inv(\w G)$ is a topological group, we have that the sequence $\{g_{n(m)}^{-1}\cdot g_{n(m+1)}\}_{m<\omega}$
converges to $p^{-1}p$, the neutral element in $G^{\rm w}$.
This takes us to a contradiction because, by Proposition \ref{compact_3},
it follows that $\{g_{n(m)}^{-1}\cdot g_{n(m+1)}\}_{m<\omega}$ must also converge to the neutral element in  $(G,\tau)$.

Therefore, the sequence $\{g_{n}\}_{n<\omega}$ is a precompact subset of $(G,\tau)$. This
implies that $p\in G$ and we are done.
\epf

\section{Locally compact groups respect other compactness-like properties}
\subsection{Countably compactness}
In this subsection we prove that LC groups respect countably compactness.

\bthm\label{CountablyCompactness1}
Let $G$ be a LC group and let $A\subseteq G$ be countably compact in $G^{\rm w}$. Then $A$ is countably compact in $G$.
\ethm
\bpf
First, we assume that $G$ is a $\sigma$-compact LC group.
If $K:=\overline A^{\, G^{\, \rm w}}$ is weakly compact, by Hughes \cite{Hughes1973}, $K$ must also be compact in the LC topology. Therefore,
the identity maps is a homeomorphism of $A^{\rm w}$ onto $A$. Therefore, $A$ is countably compact in the LC topology.
On the other hand, if $\overline A^{\, G^{\rm w}}$ is not compact, there must be some point $p\in \overline A^{\, \w G}\setminus G$.
Since the weak topology on $G$ is weaker than the original LC topology of $G$, we have that $G^{\rm w}$ is also $\sigma$-compact. That is,
there is a collection $\{K_n : n\in\omega\}$ of weakly compact subsets such that $G=\bigcup\limits_{n\in\omega}K_n$.
For every $n\in\omega$, take $f_n\in C(\w G)$ such that $0\leq f_n(x)\leq 1/2^n$ for all $x\in \w G$, $f_n(p)=0$ and $f_n(K_n)=\{1/2^n\}$.
Set $f=\sum\limits_{n\in\omega} f_n$. Then $p\in Z(f)$ and $Z(f)\cap G=\emptyset$.
Thus, the map $1/f$ is weakly continuous and, since $p\in \overline A^{\, \w G}\setminus G$,
is not bounded on $A$. This is a contradiction since $A$ is countably compact in $G^{\rm w}$. This completes the proof when $G$ is $\sigma$-compact.

In the general case, assume that $A$ is a countably compact subset of $G^{\rm w}$.
In case $A$ is precompact in $G$, it follows that $K:=\overline A^{\, G}$ is compact in $G$ and, as a consequence, is also compact in $G^{\rm w}$.
Therefore both topologies, the weak topology and the LC-topology, coincide on $K$.
Therefore, we may repeat the same argument used in the paragraph above to conclude that $A$ is countably compact in the LC topology.
Thus we may assume that $A$ is not precompact in the LC topology. As a consequence,
there is a compact neighborhood of the identity $U$ and a sequence $(a_n)\subseteq A$ such that $(a_n)$ is $U$-discrete; that is
$a_nU\cap a_mU=\emptyset$ if $n\not= m$. Set $H$ the group generated by $U\cup (a_n)$. Clearly, $H$ is a $\sigma$-compact open subgroup of $G$.
By \cite[Prop. 3.12]{FerHerTar_JFA} (cf. \cite[Lem. 3.2]{Bichteler1969}), it follows that $H$ is closed in $G^{\rm w}$ and the weak topology in $H^{\rm w}$
coincides canonically with the weak topology that $H$ inherits from $G^{\rm w}$. As a consequence, it follows that
$A\cap H$ is countably compact in $H^{\rm w}$. Since the group $H$ is $\sigma$-compact, it follows that $A\cap H$ is countably compact in $H$.
This is a contradiction because $(a_n)\subseteq A\cap H$ is $U$-discrete. This completes the proof.
\epf

\subsection{Functional boundedness}

A subset $A$ of a topological space $X$ is said to be \emph{functionally bounded} when $f_{|A}$
is bounded for every $f\in C(X)$. The topological space $X$ is a $\mu$-space when every
functionally bounded subset of $X$ is relatively compact. In this subsection, we prove that $G^{\rm w}$
is a $\mu$-space for all LC group $G$.

\blem\label{HM1}
Let $A$ be functionally bounded subset of $G^{\rm w}$ for an LC group $G$. If $p\in \overline{A}^{\, \w G}$, then $p$ belongs
to the $G_\gd$-closure of $G$ in $\w G$ and $\delta_p$ is $t_p(G)$-continuous on each
countable subset $F\subseteq I(G)$.
\elem
\bpf
Remark that $Z(f)\cap G\not=\emptyset$ for all $f\in C(\w G)$ such that $p\in Z(f)$. Indeed, if $Z(f)\cap G=\emptyset$, then $1/f$ would be weakly continuous
and, since $p\in \overline A^{\, \w G}$, we would have that $1/f$ ought to be unbounded on $A$, a contradiction.
In general, given an arbitrary $G_\delta$-open subset $N$ in $\w G$ containing $p$, it is readily seen that
there is a zero-set $Z(f)\subseteq \w G$ such that $p\in Z(f)\subseteq N$, which implies that $N\cap G\not=\emptyset$.
This verifies that $p$ belongs to $G_\gd$-closure of $G$ in $\w G$.

Now, given an arbitrary countable subset $F$ of $I(G)$, set $N:=\cap\{\varphi^{-1}(\varphi(p)) : \varphi\in F\}$.
Then $N$ is a $G_\delta$ open subset in $\w G$ and $p\in N$. As a consequence $N\cap G\not=\emptyset$.
Take $g\in N\cap G$. Then $\varphi(p)=\varphi(g)$ for all $\varphi\in F$, which yields
the continuity of $\gd_p$ on $(F,t_p(G))$.
\epf
\mkp



\bprp\label{HM2}
Let $G$ be an LC group and let $A$ be a countable subset of $G$. If $p\in \overline{A}^{\, \w G}$ and $\gd_p$ is $t_p(G)$-continuous
on each countable subset of $I(G)$, then $\gd_p$ is continuous on every compact subset of $I(G)$.
\eprp
\bpf
We define the following equivalence relation on $I(G)$: $\phi_1\sim \phi_2$ if and only if $\phi_{1|A}=\phi_{2|A}$.
Let $K$ be a compact subset of $I(G)$.
Take the \emph{quotient} map $$\pi\colon (I(G),t_p(G))\to (\, \frac{I(G)}{\sim}, t_p(A))$$ that is clearly continuous.
This means that $\pi(K)$ is compact in $(\, \frac{I(G)}{\sim}, t_p(A))$. Therefore
$\hbox{density}(\pi(K))\leq \hbox{weight}(\pi(K))=|A|=\omega$. Let $D$ be a countable subset of $K$ such that $\pi(D)$ is dense in $\pi(K)$.
Since $D$ is countable, and $p$ belongs to the $G_\gd$-closure of $G$ in $\w G$, there is $g\in G$ such that $\phi(p)=\phi(g)$ for all $\phi\in D$.
Furthermore, using the continuity of $\gd_g$ and $\gd_p$ on each countable subset of $I(G)$, it follows that we can extend $\gd_p=\gd_g$ to a continuous
map on $\overline D^{\, I(G)}$ (indeed, $\gd_p$ is continuous on $D\cup \{\phi\}$ for all
$g\in \overline D^{\, I(G)}$ and by \cite[I.57.5]{Bourbaki1971}, this implies the continuity of $\gd_p$ throughout $\overline D^{\, I(G)}$.

Consider the following diagram

\[
\xymatrix{ I(G) \ar@{>}[rr]^{\pi} \ar[dr]_{\gd_p} & & \frac{I(G)}{\sim}\ar[dl]^{\bar{\gd_p}} \\ & \mathbb D &}
\]
\mkp

where $\bar{\gd_p}$ is defined by $\bar{\gd_p}(\pi(\phi))=\gd_p(\phi)$.

Remark that $\bar{\gd_p}$ is properly defined because ${\gd_p}(\phi_1)={\gd_p}(\phi_2)$ whenever $\phi_1\sim\phi_2$.
Furthermore, it holds that $\bar{\gd_p}$ is continuous on $\pi(\overline D^{\, I(G)})$. Indeed, in order to verify this, it will suffice to prove that
for every $\phi\in \overline D^{\, I(G)}$ and each net $(\phi_i)$ in $\overline D^{I(G)}$  such that $\pi((\phi_i))$ converges to  $\pi(\phi)$,
there is a subnet $(\phi_m)$ such that $\gd_p((\phi_m))$ converges to $\gd_p(\phi)$.

Now, since $\overline D^{\, I(G)}$ is compact, it follows that there is a subnet $(\phi_m)$ converging to $\phi'\in \overline D^{\, I(G)}$. By the continuity of $\pi$
the net $\pi((\phi_m))$ must converge to $\pi(\phi')$ but, by our previous assumption, will also converge to $\pi(\phi)$.
This means that $\phi\sim\phi'$ and, as a consequence, $\gd_p(\phi)=\gd_p(\phi')$. The continuity of $\gd_p$ on $\overline D^{\, I(G)}$ implies
the convergence of $\gd_p((\phi_m))$ to $\gd_p(\phi)=\gd_p(\phi')$. Bearing in mind the definition of $\bar\gd_p$, the continuity
of this map has been proved.

Now, since $\pi(D)$ is dense in $\pi(K)$, it follows that $\pi(\overline D^{\, I(G)})=\pi(K)$. Thus, we have proved the continuity of $\bar\gd_p$ on $\pi(K)$.
The commutativity of the diagram

\[
\xymatrix{ K \ar@{>}[rr]^{\pi} \ar[dr]_{\gd_p} & & \pi(K)\ar[dl]^{\bar{\gd_p}} \\ & \mathbb D &}
\]
\mkp

implies the continuity of $\gd_p$ on $K$, which completes the proof.
\epf
\mkp

Now follows the main result of this paper. It extends to non-necessarily abelian groups a previous
result of Trigos-Arrieta [23] for locally compact Abelian groups.

\bthm\label{FunctionalBoundedness}
Let $G$ be a LC group. Then the group $G^{\rm w}$ is a $\mu$-space.
\ethm
\bpf
We must verify that every closed functionally bounded subset of $G^{\rm w}$ is compact in $G^{\rm w}$.
Let $A$ be a closed functionally bounded subset of $G^{\rm w}$. If $A$ is countably compact in $G^{\rm w}$
then, by Theorem \ref{CountablyCompactness1}, it follows that $A$ is closed and countably compact in $G$.
Hence, we have that $A$ is compact in $G$ and, as a consequence, in $G^{\rm w}$.

Therefore, from here on, we  assume that $A$ is not countably compact in $G^{\rm w}$ without loss of generality.
This implies that there is some sequence $(a_n)\subseteq A$ that has no closure points in $G^{\rm w}$.
Let $p\in \overline{(a_n)}^{\, \w G}\setminus G$. By Lemma \ref{HM1}, $p$ belongs to the $G_\gd$-closure of $G$
in $\w G$ and $\gd_p$ is $t_{(G)}$-continuous on each countable subset
of $I(G)$. Then, by Proposition \ref{HM2}, we deduce that $\gd_p$ is $t_p(G)$-continuous on each compact subset of $I(G)$.

\noindent The proof now requires a case-study approach on the structure of the group $G$. 

\begin{enumerate}[(1)]
\item $G$ is $\sigma$-compact and metrizable and, therefore, a Polish LC group. In this case, the space
$(I(G),t_k(G))$ is metrizable and, as a consequence, a $k$-space. Since we have verified that
$\gd_p$ is $t_p(G)$-continuous on each compact subset of $I(G)$, it follows that
$\gd_p$ is continuous on $I(G)$. Applying Akenmann duality thorem \cite{Akemann1972}, it follows that
$p\in G$ and, as a consequence, that $\overline A^{\, \w G}\subseteq G$.\medskip

\item Suppose first that $(a_n)$ is not precompact in $G$. Then we may assume, with some notational abuse,
that $(a_n)$ is $U$-discrete for some compact neighborhood of the identity $U$. Therefore,
by Kakutani-Kodaira's theorem, there exists a normal, compact $K$ of $G$ such that $K\subseteq U$ and
$G/K$ is metrizable, and consequently Polish.
Let $p:G\rightarrow G/K$ be the quotient map, which is also continuous for the the weak topologies
$p:G^{\rm w}\rightarrow (G/K)^{\rm w}$, since every pure positive definite map on $G/K$ can be lifted canonically to $G$.
Therefore $p((a_n))$ is a functionally bounded subset in $(G/K)^{\rm w}$. By (1) $p((a_n))$ is
relatively compact in $G/K$, which is a contradiction because $(a_n)$ was assumed to be $U$-discrete
and $K\subseteq U$. Thus, we may assume, without loss of generality, that $(a_n)$ is precompact in $G$.
This means that $(a_n)$ is relatively compact in $A$, which is a closed subset of $G$.
Therefore, we have again that $\overline A^{\, \w G} \subseteq G$.

\item $G$ is LC. Again, we suppose first that $(a_n)$ is not precompact in $G$ and, with some notational abuse,
that $(a_n)$ is $U$-discrete for some compact neighborhood of the identity $U$.
Consider the subgroup $H$ generated by $U\cup\lbrace g_n\rbrace_{n<\omega}$, which $\sigma$-compact and open in $G$.
Since $\overline H^{\, \w G}$ is canonically homeomorphic to $\w H$, we may identify
$\overline{\lbrace g_n\rbrace_{n<\omega}}^{\, \w H}$ with $\overline{\lbrace g_n\rbrace_{n<\omega}}^{\, \w G}$. Hence
we may assume, without loss of generality, that $G$ is $\sigma$-compact and the proof follows from (2).
The case $(a_n)$ is precompact follows also from (2).
\end{enumerate}

Thus, we have proved that $\overline A^{\, \w G}\subseteq G$ for evey functionally bounded subset $A$ of $G^{\rm w}$, which proves that
$G^{\rm w}$ is a $\mu$-space.
\epf

\bcor
Every LC group respects pseudocompactness.
\ecor


\end{document}